\documentclass[12pt]{article}

\usepackage{amsmath, amsbsy, amsthm, amsfonts, array,graphics}

\def\XXint#1#2#3{{\setbox0=\hbox{$#1{#2#3}{\int}$}
\vcenter{\hbox{$#2#3$}}\kern-.5\wd0}}

 \DeclareMathOperator{\der}{d}
 \DeclareMathOperator{\Vol}{Vol}
 
 \DeclareMathOperator{\B}{B}
 
 \DeclareMathOperator{\Ric}{Ric}
 \DeclareMathOperator{\Rm}{Rm}

 \newtheorem{thm}{Theorem}[section]
 \newtheorem{cor}[thm]{Corollary}
 \newtheorem{lem}[thm]{Lemma}

 \newtheorem{rem}[thm]{Remark}
 \numberwithin{equation}{section}

\title
{\textbf{ \normalsize{Smoothing Riemannian Metrics with Bounded Ricci Curvatures in
Dimension Four, II}}}
\author
{\normalsize{Ye Li}\\
\it{\footnotesize{Mathematics Department, Duke University, Durham, NC 27708}} \\
\it{\footnotesize{Email address: yeli@math.duke.edu}}
}
\date{}

\begin{document}
\maketitle
\begin{center}
\begin{minipage}{12cm}
    \hspace{0.5cm}\textbf{Abstract:} This note is a continuation of the author's paper \cite{Li}. We prove that if
    the metric $g$ of a 4-manifold has bounded Ricci curvature and the curvature has
    no local concentration everywhere, then it can be smoothed to a metric with bounded sectional curvature. Here we don't assume the bound for local Sobolev constant of $g$ and hence this smoothing result can be applied to the collapsing case.
\end{minipage}
\end{center}

\vskip 0.2cm
\begin{center}
\begin{minipage}{12cm}

\end{minipage}

\end{center}
\vskip 0.2cm

\section{Introduction}

It is well known that if a Riemannian manifold has bounded sectional
curvature, then at least its local structure is better understood
than that with weaker curvature bounds. A natural question is how to generalize the results which hold for spaces with bounded curvatures to those with bounded Ricci curvatures. Therefore, it is important to discover whether we can deform
or smooth a metric with bounded Ricci curvature to a metric with bounded
curvature. In particular, if one wants to study the compactness theory for a sequence of metrics  with uniform bounded Ricci curvatues, the first step is to construct   ``good'' coordinates of definite size, say harmonic coordinates, at the place where the metric has no local curvature concentration such that the space will not develop singularities. However, if, in addition,  we don't assume the noncollapsing condition, such coordinates will not exist. A standard argument dealing with the collapsing situation is to lift the metric to the tangent space and construct ``good'' coordinates there, but, we cannot immediately do so under the assumption of bounded Ricci curvature. Therefore, it is necessary to smooth the metric to the one with bounded curvature.

In this note we study the problem whether we
can smooth a metric with bounded Ricci curvature to that with
bounded curvature in dimension 4. Due to the reason mentioned above, our result does not depend on the noncollapsing condition.

Let $X=S^2\times S^2 / \sigma$, where $\sigma$ denotes the Cartan
involution that rotates $\pi$ on each sphere. We can resolve the
four resulting singularities to obtain
$M=\mathbb{C}P^2\#5\overline{\mathbb{C}P^2}$ which admits a sequence
of K\"{a}hler-Einstein metrics converging to $X$. With regard to
this example, in general, we cannot smooth a metric with bounded
Ricci curvature  to the one with bounded curvature.
Consequently, additional geometric conditions must be needed.

To the author's knowlege, two methods are known in the literature to do the smoothing.
One is the heat flow method. If we assume that the initial metric 
has bounded curvature, then we can show the short-time existence of 
Ricci flow and obtain the covariant derivative bounds for the 
curvature tensors along the Ricci flow,  for example, \cite{BMR}, \cite{Ba} and
\cite{Sh}. In the case that the inital metric has only bounded Ricci curvature, some additional condition is needed for the smoothing procedure. 
In \cite{DWY}, Dai, Wei and Ye studied how to do smoothing on a compact
manifold with bounded Ricci curvature  and conjugate radius bound. Here the conjugate radius bound ensures that one can lift the metric to the tangent space and hence obtain certain initial conditions to run the flow. Thus to smooth the metrics one needs to find out the suitable conditions on the geometry of the space to rule out the possible singularities developing, which is not convenient in general. So one may consider the smoothing procedure locally, which can be easily performed in practice.  In \cite{Ya}, Yang introducted the local Ricci flow and proved its short-time existence under the assmptions that the initial metric has some suitable integral curvature bounds and satisfies noncollapsing condition, and as a consequence, the metric can be smoothed in a given geodesic ball with a definite volume. Since Yang's result is purely local, we don't need the global geometric assumtions in the space.

The other way is the embedding method, for
example, \cite{CG}, \cite{Ab} and  \cite{PWY}. The basic idea of the embedding method is that if the injectivity radius of $M$ has a positive lower bound, we can embed $M$ into $L^2(M)$ by a map $I$ involving the distance function such that a sufficiently smooth submanifold can be found in a neighbor of $I(M)$. We can perform this construction in each coordinate chart of $M$ and glue them together to get the desired smooth metric. If we don't have the assumption on the injectivity radius, then each small neighborhood of $M$ is a quotient of a Riemannian manifold satisfying injectivity radius lower bound by an action of a pseudofundamental group.
However, embedding method may not
be as convenient as the heat flow method if the metric has no bounded
sectional curvature. 

Recently, in \cite{Xu}, using local Ricci flow developed by Yang \cite{Ya}, Xu proved the short-time existence of Ricci flow
under the assumptions that the initial metric satisfies the volume
doubling property, local Sobolev constant bound and some integral
bounds on curvature such that the curvature has local concentration nowhere and thus the space will not develop singularities. The purpose of this note is two-fold. Firstly, if the curvature has
no local concentration everywhere, one can prove the short-time
existence of Ricci flow directly without using the local Ricci flow.
As a matter of face, we can employ a covering argument as in \cite{DWY} to obtain
the corresponding energy estimates needed in the short-time
existence of Ricci flow and thus avoid using the local Ricci flow.
This can then simplify the proof given in \cite{Xu}. Secondly, based on the recent breakthrough by Cheeger and Tian in the work of collasping Einstein 4-manifolds, \cite{ChGi}, if the initial metric has bounded Ricci
curvature, we can remove the assumption on the local Sobolev
constant in dimension 4, or equivalently, the lower volume bound of
each geodesic ball involved, and thus the result can be applied to the collapsing case in dimension 4. Since the proof in \cite{ChGi} depends on Gauss-Bonnet-Chern formula, it is not known if one can generalize the result to higher dimensions.

For convenience, we will give a brief description of removing the noncollapsing assumption in dimension 4.  In \cite{ChGi},  Cheeger and Tian obtained the following technical result: Given a 4-dimensional Einstein manifold with Einstein constant $|\lambda|\le3$, if there exists some $\varepsilon>0$ such that $\int_{B_r(x)}|\Rm|^2<\varepsilon$, then
for some constant $c>0$,
$$\frac{\Vol(B_{cr}(\underline{x}))}{\Vol(B_{cr}(x))}\int_{B_{cr}(x)}|\Rm|^2< \varepsilon,$$
where $\underline{x}$ is a point in the simply connected space of
constant curvature $-1$. In \cite{Li}, the author generalized it to the 4-dimensional bounded Ricci curvature case. First we need to smooth the metric in some suitable local scale to obtain
the curvature bound for some nearby metric. Here the estimates are not uniform and depend on the chosen local scale. Then by using a local equivariant version of good chopping and an iteration technique one can show the following key estimate $$\frac{\Vol(B_{r}(\underline{x}))}{\Vol(B_{r}(x))}\int_{B_{r}(x)}|\Rm|^2\le C,$$
where $C$ is a definite constant. This result is nontrivial in the collapsing case, since the volume of a geodesic ball is arbitrarily small. If we can lift the metric to the tangent space, the corrsponding $L^2$-norm of curvature may not be small and even worse it may be unbounded. However the above estimate indicates that although the $L^2$-norm of curvature may not be small, it is still bounded and thus rule out the worse situation. Finally, one can show that the 
quantity $$\frac{\Vol(B_{r}(\underline{x}))}{\Vol(B_{r}(x))}\int_{B_{r}(x)}|\Rm|^2$$ can be sufficiently small if one shrinks the ball to a smaller concentric one, whose radius is comparable to $r$. The proof
used the Gauss-Bonnet-Chern formula, an estimate on the
transgression form in terms of volume growth rate, and a controlled and smooth approximation of the
distance function.

According to \cite{Bu}, \cite{Va}, \cite{Sa1}, \cite{Sa2} and \cite{An}, etc., there exists a constant
$C$ depending only on the bounds of Ricci curvature and the
dimension of $M$ such that the Sobolev constant for the geodesic
ball $B_r(x)$ can be controlled as follows: for $r\le 1$,
$$C_s(B_r(x))\le C\left(\frac{r^4}{\Vol(B_r(x))}\right)^{\frac {1}{2}}.$$
This implies that the $L^2$-norm of curvature actually can be made arbitrarily small against local Sobolev constant.

Our main result is the following

\begin{thm}
Let $(M,g_0)$ be a complete noncompact Riemannian 4-manifold. There
exist constant $\varepsilon$ and $C_1$ such that for $r\le 1$, if
$$\int_{\B_r(x)}|\Rm(g_0)|^2\der V_{g_0}\le \varepsilon,\hbox{ for any }x\in M,$$
and
$$|\Ric(g_0)|\le K,$$
then the Ricci flow
\begin{align*}
    \left\{
    \begin{array}{ll}
       \dfrac{\partial g}{\partial t}&=-2\Ric(g),\\
       g(0)&=g_0.\\
    \end{array}
    \right.
\end{align*}
has a smooth solution for $t\in[0,T)$, where
$$T\ge C_1\cdot\min \left(r^2, K^{-1} \right).$$
Moreover, for $t\in(0,T)$, the Riemannian curvature tensor satisfies
the following bound,
\begin{align*}
\|\Rm\|_{\infty}\le C_2t^{-1}.
\end{align*}
Here $\varepsilon$, $C_1$ and $C_2$ only depend on the dimension of
$M$.
\end{thm}

\vskip 0.2cm

This note is organized as follows. In Section 2, we show the Moser's iteration for linear heat equations. In Section 3, we will study the short-time existence of nonlinear equaiton and apply the result to Ricci flow.

\vskip 0.2cm

\vskip 0.2cm

\noindent \textbf{Acknowledgement} The author would like to thank
Professors Jeff Cheeger and Gang Tian for many helpful discussions on their work \cite{ChGi} and their support. The author is also grateful to Professor
Laurent Saloff-Coste for discussion on the literature of the local Sobolev constant bounds.

\vskip 1cm

\section{Moser's Iteration for Linear Heat Equations}

In this section we study Moser's weak maximum principle for linear
equations. We mainly follow the lines in \cite{Ya}. The difference is that the equation discussed here is not local. 

Fix a geodesic ball $B_r(x)\subset M^4$ and a smooth compactly
supported function $\phi\in C_0^{\infty}(B_r(x))$.

Let $g(t),0\le t \le T$, be a 1-parameter family of smooth
Riemannian metrics. Let $\nabla$ denote the covariant
differentiation with respect to the metric $g(t)$ and $-\triangle$
be the corresponding Laplace-Beltrami operator. Let $A>0$ be a
constant that satisfies the standard Sobolev inequality
\begin{align}\label{5}
\left(\int_{B_r(x)}f^{4}\der V_g\right)^{\frac 12}\le A
\int_{B_r(x)}|\nabla f|^2\der V_g,\ f\in C_0^{\infty}(B_r(x)),
\end{align}
with respect to each metric $g(t), 0\le t\le T$.

Assume that for each $t\in [0,T]$,
$$\frac{1}{2}g_{ij}(0)\le g_{ij}(t)\le 2 g_{ij}(0)\ \ \text{on}\ B_r(x).$$
All geodesic balls in this section are defined with respect to the
metric $g(0)$, and therefore, are fixed open subsets of $M$, and
independent of $t$.

We want to study the heat equation:
\begin{align}\label{1}
\frac{\partial f}{\partial t}\le \triangle f+uf,\ 0\le t \le T,
\end{align}
where $f$ and $u$ are nonnegative functions on $B_r(x)\times [0,T]$,
such that
\begin{align}\label{2}
\frac{\partial}{\partial t}\der V_g\le c\cdot u\der V_g
\end{align}
and
\begin{align}\label{3}
\left(\int_{B_r(x)}u^3\right)^{\frac 13}\le \mu t^{-\frac 13}.
\end{align}

\vskip 0.2cm

The following two lemmas are the consequences of direct
computations.
 \vskip 0.2cm
\begin{lem}\label{lem1}
Given $p>1,\  \phi\in C_0^{\infty}(B_r(x)),\ f\in C^{\infty}(M),\
f\ge 0$,
\begin{align*}\int_{B_r(x)}|\nabla(\phi f^{\frac{p}{2}})|^2\le \frac{p^2}{2(p-1)}&\int_{B_r(x)}\phi^2f^{p-1}(-\triangle f)\der V_g\\
+&\left(1+\frac{1}{(p-1)^2}\right)\int_{B_r(x)}|\nabla \phi|^2f^p\der
V_g.
\end{align*}
\end{lem}

\vskip 0.2cm

 \noindent\textit{Proof:}  Using integration by parts,
we have
\begin{align*}
\int |\nabla (\phi f^{\frac{p}{2}})|^2 =& -\int \phi f^{\frac p2} \triangle (\phi f^{\frac p2})\\
    =&\frac p2 \int \phi^2 f^{p-1} (-\triangle f)+ \int f^p|\nabla \phi|^2-\frac{p(p-2)}{4}\int \phi^2 f^{p-2}|\nabla f|^2\\
    =&\frac {p^2}{2(p-1)} \int \phi^2 f^{p-1}(-\triangle f)+\frac{p}{2(p-1)}\int \phi^2 f^{p-1} \triangle f \\
     & +\int f^p |\nabla \phi|^2-\frac{p(p-2)}{4}\int \phi^2 f^{p-2} |\nabla f|^2.
\end{align*}
On the other hand, by Cauchy inequality,
\begin{align*}
\frac{p}{2(p-1)}\int \phi^2 f^{p-1}\triangle f &= -\frac{p}{2(p-1)}\int \nabla (\phi^2 f^{p-1})\nabla f\\
    &=-\frac{p}{p-1}\int \phi f^{p-1} \nabla \phi \nabla f-\frac{p}{2}\int \phi^2 f^{p-2}|\nabla f|^2\\
    &\le \frac{1}{(p-1)^2}\int f^p |\nabla \phi|^2+\frac{p^2}{4}\int \phi^2f^{p-2}|\nabla f|^2-\frac p2 \int \phi^2 f^{p-2}|\nabla f|^2\\
    &=\frac{1}{(p-1)^2}\int f^p|\nabla \phi|^2+\frac{p(p-2)}{4}\int \phi^2 f^{p-2}|\nabla f|^2.
\end{align*}
This proves the lemma. $\Box$

\vskip 0.2cm

\begin{lem}\label{lem2}
Suppose that $f$ and $u$ are nonnegative functions on
$B_r(x)\times[0,T)$ which satisfy (\ref{1}), (\ref{2}) and
(\ref{3}). For $p>1,$ we have
\begin{align}\label{4}
\frac{\partial}{\partial t}\int \phi^2f^p+\frac{p-1}{p}\int
|\nabla(\phi f^{\frac{p}{2}})|^2 \le C_p\int
|\nabla\phi|^2f^p+C_p\mu^3A^2t^{-1}\int \phi^2f^p.
\end{align}
\end{lem}

\vskip 0.2cm

\noindent\textit{Proof:} By Lemma \ref{lem1}, we have
$$
\frac{\partial}{\partial t}\int \phi^2f^p+2\left(1-\frac 1p\right)\int |\nabla(\phi f^{\frac p2})|^2\\
    \le C_p\int |\nabla \phi|^2f^p+(p+c)\int u\phi^2f^p.
$$
Using H$\ddot{o}$lder, Sobolev, Cauchy inequalities, and (\ref{3}),
we see that
\begin{align*}
\int u\phi^2f^p\le &\left(\int  u^3\right)^{\frac
13}\left(\int(\phi^2f^p)\right)^{\frac 13}
                \left(\int \phi^4f^{2p}\right)^{\frac 13}\\
    \le & \mu t^{-\frac 13}\left(\int(\phi^2f^p)\right)^{\frac 13}\cdot A^{\frac 23} \left(|\nabla(\phi f^{\frac p2})|^2\right)^{\frac 23}\\
    \le & (\mu t^{-\frac 13})^3\varepsilon^{-\frac 13}\int \phi^2f^p+\varepsilon^{\frac 23}A^2\int |\nabla(\phi f^{\frac p2})|^2.
\end{align*}

Choosing $\varepsilon$ so that $\varepsilon^{\frac
23}A^2=\frac{p-1}{p}$ , we complete the proof of lemma \ref{lem2}.
$\Box$

\vskip 0.2cm

Now given $0<\tau<\tau '<T$, let
\begin{align*}
    \psi(t)=\left\{
    \begin{array}{cc}
       0, &0\le t\le \tau,\\
       \dfrac{t-\tau}{\tau '-\tau}, &\tau\le t\le \tau ',\\
       1, &\tau '\le t\le T.
    \end{array}
    \right.
\end{align*}
Multiplying (\ref{4}) by $\psi$, we have
\begin{align*}
\frac{\partial}{\partial t}\left(\psi \int \phi^2f^p\right)+&\psi\int |\nabla (\phi f^{\frac p2})|^2\\
    \le& C_p\psi\int|\nabla\phi|^2f^p+\left(\hat{C}(t)\psi+|\psi '| \ \right)\int \phi^2f^p,
\end{align*}
where $\hat{C}(t)=C_p\mu^3 A^2 t^{-1}$. Integrating this with
respect to $t$ from $\tau$ to $t$ and throwing away $\int_{\tau}^{\tau'}\psi|\nabla (\phi f^{ \frac
p2})|^2$ on the left-hand side, we obtain

\vskip 0.2cm
\begin{lem}\label{lem3}
For $\tau '\le t\le T$, we have
\begin{align*}
\int_t \phi^2f^p+\int_{\tau '}^{t}\int |\nabla (\phi f^{ \frac
p2})|^2 \le
C_p\int_{\tau}^T\int|\nabla\phi|^2f^p+\left(\hat{C}(\tau')
 +\frac{1}{\tau '-\tau}\right)\int_{\tau}^T\int \phi^2f^p.
\end{align*}
\end{lem}

\vskip 0.2cm

Now given $p>1,\ 0\le \tau <T$, denote
$$H(p,\tau,r)=\int_{\tau}^{T}\int_{B_r(x)}f^p.$$

\vskip 0.2cm
\begin{lem}\label{lem4}
Given $p\ge p_0$, $0\le \tau < \tau ' < T$ and $r'<r$
\begin{align*}
H\left(\frac 32 p,\tau',r'\right)\le A\left(\hat{C}(\tau')
 +\frac{1}{\tau '-\tau}+\frac{C_1}{(r-r')^2}\right)^{\frac 32}
 H(p,\tau,r)^{\frac
32}.
\end{align*}
\end{lem}

\noindent\textit{Proof:} Choosing a suitable cut-off function $\phi$
and noticing $|\nabla\phi|_t\le 2|\nabla\phi|_0$, we have
\begin{align*}
H\left(\frac 32 p,\tau',r'\right)\le &\int_{\tau '}^{T}\int_{B_r(x)} (\phi^2f^p)^{\frac 32}\\
        \le& \int_{\tau '}^{T} \left(\int_{B_r(x)} \phi^2f^p\right)^{\frac 12}\left(\int_{B_r(x)} \phi^4f^{2p}\right)^{\frac 12}\der t\\
        \le& \left(\sup_{\tau '\le t\le T}\int_{B_r(x)} \phi^2f^p\right)^{\frac 12}\ A\ \int_{\tau '}^T\int_{B_r(x)} |\nabla (\phi f^{\frac p2})|^2 \der
        t\\
        \le& A\left(4\int_{\tau}^T\int|\nabla\phi|^2f^p+\left(\hat{C}(\tau')
 +\frac{1}{\tau '-\tau}\right)\int_{\tau}^T\int
 \phi^2f^p\right)^{\frac 32}\\
 \le& A\left(\hat{C}(\tau')
 +\frac{1}{\tau '-\tau}+\frac{C_1}{(r-r')^2}\right)^{\frac 32}\left(\int_{\tau}^T\int \phi^2f^p\right)^{\frac 32}
\end{align*}
This proves the lemma. $\Box$

The following theorem is the consequence of Moser's iteration.

 \vskip 0.2cm

\begin{thm}\label{thm1}
Let $f$ and $u$ be non-negative functions on $B_r(x)\times[0,T)$
satisfying (\ref{1}), (\ref{2}) and (\ref{3}). Then for $t\in[0,T)$
and $p_0>2$,
$$|f(x,t)|\le CA^{\frac {2}{p_0}}\left((1+A^2\mu^3)t^{-1}+\frac{1}{r^2}\right)^{\frac{3}{p_0}}
    \left(\int_0^T\int_{B_r(x)} f^{p_0}\right)^{\frac{1}{p_0}},$$
where $C$ depends on the dimension of $M$, $p_0$.
\end{thm}

\vskip 0.2cm

\noindent \textit{Proof: } Denote $\nu=\frac 32$. Fix $0<t<T$, and
set
\begin{align*}
p_k=&\ p_0\nu^k,\\
\tau_k=&\ t(1-\nu^{-k-1}),\\
r_k=&\frac r2\left(1+\nu^{-k/2}\right),\\
\Phi_k=&\ H(p_k,\tau_k,r_k)^{\frac {1}{p_k}}.
\end{align*}
Applying Lemma \ref{lem4},
\begin{align*}
H(p_{k+1},\tau_{k+1},r_{k+1})\le AC
\left((1+A^2\mu^3)t^{-1}+\frac{1}{r^2}\right)^{\nu}\nu^{k\nu}
    H(p_k,\tau_k,r_k)^{\nu}.
\end{align*}
Therefore,
$$\Phi_{k+1}\le (AC)^{\frac{1}{p_{k+1}}}\left((1+A^2\mu^3)t^{-1}+\frac{1}{r^2}\right)^{\frac{1}{p_k}}
   \nu^{\frac{k}{p_k}}\Phi_k,$$
Hence,
$$\Phi_{k+1}\le (AC)^{\frac{\sigma_{k+1}}{p_0}}\left((1+A^2\mu^3)t^{-1}+\frac{1}{r^2}\right)^{\frac{\sigma_k}{p_0}}
   \nu^{\frac{\sigma_{k}^{\prime}}{p_0}}H(p_0,0,r)^{\frac {1}{p_0}},$$
where $\sigma_k=\sum_{i=0}^{k}\nu^{-i},\
\sigma^{\prime}_k=\sum_{i=0}^{k}i\nu^{-i}.$ Letting $k\rightarrow
\infty$, we obtain
$$|f(x,t)|\le CA^{\frac {2}{p_0}}\left((1+A^2\mu^3)t^{-1}+\frac{1}{r^2}\right)^{\frac{3}{p_0}}
    \left(\int_0^T\int_{B_r(x)} f^{p_0}\right)^{\frac{1}{p_0}}.$$
Now let $T\rightarrow t$.  This proves the theorem.  $\Box$

\vskip 1cm

\section{Short-time Existence for Ricci Flow}

In this section we study the short-time existence of Ricci flow and hence obtain the smoothing result. We will follow the lines in \cite{Ya}, together with the covering argument in \cite{DWY} to get the desired energy estimates in Ricci flow.

Let $M$ be a complete noncompact manifold with Riemannian metric
$g_0$. Consider the following evolution equation
\begin{align}\label{8}
    \left\{
    \begin{array}{ll}
       \dfrac{\partial g}{\partial t}&=-2\Ric(g),\\
       g(0)&=g_0.\\
    \end{array}
    \right.
\end{align}

It is easy to check that the curvature tensor $\Rm$ and Ricci tensor
$\Ric$ satisfy the following equations respectively,
\begin{align}\label{9}
\frac{\partial \Rm}{\partial t}=\triangle
\Rm+Q_1(\Rm,\Rm),
\end{align}
and
\begin{align}\label{10}
\frac{\partial \Ric}{\partial t}=\triangle
\Ric+Q_2(\Rm,\Ric),
\end{align}
where $Q_i$ are multi-linear functions of
their arguments, $i=1,2$. Their definitions depend only on the
dimension of $M$.

\vskip 0.2cm

\begin{thm}\label{thm4}
There exist constant $C_1$ and $C_2$ such that if for all $x\in M$
$$\left(\int_{B_r(x)}|\Rm(g_0)|^2\der V_{g_0}\right)^{\frac 12}\le [\ C_1C_s(B_r(x))\ ]^{-1}$$
and
$$|\Ric(g_0)|\le K,$$
then the equation (\ref{8}) has a smooth solution for $t\in[0,T)$,
where
$$T\ge C_2\cdot\min \left(r^2, K^{-1} \right).$$
Moreover, for $t\in(0,T)$, the Riemannian curvature tensor satisfies
the following bound,
\begin{align}\label{11}
\|\Rm\|_{\infty}\le \frac{C_3}{t}.
\end{align}
Here $C_1$, $C_2$ and $C_3$ only depend on the dimension of $M$.
\end{thm}

\vskip 0.2cm

\begin{rem}
By the same argument as in the proof, we can show that Theorem \ref{thm4} holds for $n\ge 3$, where the $L^2$-norm of the curvature should be replaced by $L^{\frac n2}$-norm.
\end{rem}

\vskip 0.2cm
\begin{rem}
The assumption of the Ricci curvature is only used to guarantee that we have a standard covering property on $M$. So if we assume that $(M,g_0)$ has such covering property, then we can weaken the assumption of Ricci curvature by $\int_{B_r(x)}|\Ric(g_0)|^p\le K$, for each $x\in M$, where $p>\frac n2$. Note that Corollary \ref{cor11} will give the required estimate on Ricci curvature along the Ricci flow.
\end{rem}
\vskip 0.2cm

To prove Theorem \ref{thm4}, we first show a result for the scalar
function. In the following we always assume that $(M,g_0)$ is a
Rimannian manifold with bounded Ricci curvature and that for each
$t\in[0,T]$ and $x\in M$,
$$\frac 12g_0\le g(t)\le 2g_0,$$
$$\left(\int_{B_r(x)}f^{4}\der V_g\right)^{\frac 12}\le A
\int_{B_r(x)}|\nabla f|^2\der V_g,\ f\in C_0^{\infty}(B_r(x)).$$

\begin{thm}\label{thm2}
Let $f\ge 0$ solve
\begin{align}\label{6}
\frac{\partial f}{\partial t}\le \triangle f+C_0 f^2,\ 0\le t\le T,
\end{align}
on $M\times[0,T)$.  Assume that
$$\frac{\partial}{\partial t}\der V_{g(t)}\le cf\der V_{g(t)}$$
and that
$$\left(\int_{B_r(x)}f^2_0\right)^{\frac 12}\le (6C_0 A)^{-1},\hbox{ for any }x\in M$$
where $f_0(x)=f(x,0)$.  Then
$$|f(x,t)|\le C_1t^{-1},$$
where $0<t<\min (T, C_2r^2)$, $C_1$ and $C_2$ depend on the
dimension of $M$ and $C_0$.
\end{thm}

\noindent \textit{Proof: } Let $[0,T']\subset [0,T)$ be the maximal
interval such that
$$e_0=\sup_{x\in M,\ 0\le t\le T'}\left(\int_{B_r(x)}f^2\right)^{\frac 12} \le (3C_0A)^{-1}.$$
By a direct calculation, we have, for $0\le t\le T'$,
\begin{align}\label{15}
\frac{\partial}{\partial t}\int \phi^{q+2}f^p+\int|\nabla(\phi f^{p/2})|^2\le
C_{p,q}\|\nabla \phi\|^2_{\infty}\int \phi^qf^p.
\end{align} 
Since $g_0$ has bounded Ricci curvature and all the metrics $g(t)$ are equivalent,
by a standard covering theorem there exists a definite constant $N$
such that
$$B_{2r}(x)\subset \cup_{i=1}^{N}B_r(y_i),\ y_i\in B_{\frac 32 r}(x).$$
Then we can choose $\phi$ such that the following holds
$$\frac{\partial}{\partial t}\int_{B_r(x)} f^2\le 2N\|\nabla \phi\|^2_{\infty}e_0^2.$$
Integrating this, we obtain
$$\int_{B_r(x)} f^2 \der V_g\le \int_{B_r(x)} f_0^2 \der V_{g(0)}+2N\|\nabla \phi\|^2_{\infty}e_0^2t.$$
Therefore,
$$(1-2N\|\nabla \phi\|^2_{\infty}t)e_0^2\le (6C_0 A)^{-1}$$
Since $x$ is chosen arbitrary, for $T'< \frac{3}{8N}\|\nabla
\phi\|_{\infty}^{-2}$, we have
$$e_0< (3C_0A)^{-1}.$$
This contradicts the assumed maximality of $[0,T']$.  We can
therefore assume that $T'\ge \min (C_2r^2,T)$.

Multiplying (\ref{15}) by $t$ and then integrating with respect to $t$, we have
$$\int\phi^{q+2}f^p+\int_0^t\int|\nabla(\phi f^{p/2})|^2\le Ct^{-1}\int_0^t\int \phi^qf^p$$
Then the covering argument yields that
\begin{align*}
\int_{B_r(x)}f^3 \le& \sum_{i=1}^{N}\int_{B_r(y_i)}\phi_i^4f^3\\
\le&Ct^{-1}\int_0^t\left(\int_{B_r(y_i)}f^2\right)^{\frac
12}\left(\int_{B_r(y_i)}(\phi_if)^4\right)^{\frac 12}\\
\le&Ce_0At^{-2}\int_0^t\int_{B_r(y_i)} f^2\\
\le& CAe_0^3t^{-1},
\end{align*}
where $\phi_i$ is a cutoff function in $B_r(y_i)$.

Setting $\mu^3=CNAe_0^3$, we can apply Theorem \ref{thm1}, which
still holds, as $p_0\rightarrow 2$.

On the other hand, for $t \in [0, T']$
\begin{align*}
A^2\mu^3 \le  CNA^3e_0^3
         \le CNA^3(3C_0A)^{-3}
         \le  CNC_0^{-3}.
\end{align*}
We then obtain the desired estimate.  $\Box$

\vskip 0.2cm

The argument also implies the following

\begin{cor}\label{cor11}
Let $f$ satisfy the assumptions of Theorem \ref{thm2}.  Then given
$u\ge 0$ such that
$$\frac{\partial u}{\partial t}\le \phi^2(\triangle u+c_0fu),$$
the following estimate holds for $0\le t< \min (T, C_2r^2)$,
$$|u(x,t)|\le C_1A^{\frac23}t^{-\frac23}\left(\int_{B_{3r}(x)}u^3_0\right)^{\frac13},$$
where $u_0(x)=u(x,0)$.
\end{cor}

\vskip 0.2cm

\noindent\textit{Proof: } Using the covering argument and the proof
of Lemma \ref{lem2}, we obtain

$$\frac{\partial}{\partial t}\int \phi^2u^3
+\frac23\int |\nabla(\phi u^{\frac 32})|^2 \le C_3\int
u^3+C_4A\left(\int f^2\right)^{\frac 12}\int |\nabla(\phi u^{\frac
32})|^2.$$ This implies that
$$\int_{B_r(x)}u^3\le C\int_{B_{3r}(x)} u^3_0.$$
By the above theorem, we have
$$\int_{B_r(x)}f^3\le\mu^3t^{-1}.$$
Thus the result follows from Theorem \ref{thm1} with $p_0=3$. $\Box$

\vskip 0.2cm

\vskip 0.2cm

\noindent\textit{Proof of Theorem \ref{thm4}: }It is well known that
the equation (\ref{8}) has a smooth solution on a sufficiently small
time interval starting at $t=0$. Let $[0,T_{\max})$ be a maximal
time interval on which (\ref{8}) has a smooth solution and such that
the following hold for each metric $g(t)$ and $x\in M$
\begin{align}
\label{12} \|f\|^2_4\le 4 A_0 \|\nabla f\|^2_2,\ f\in
C_0^{\infty}(B_r(x)),
\end{align}
\begin{align}
\label{13} \frac 12 g_0 \le g(t) \le 2 g_0,
\end{align}
\begin{align}
\label{14} \left(\int_{B_r(x)}|\Rm|\right)^{\frac 12}\le
2(C_1A_0)^{-1},
\end{align}

Suppose that $T_{\max}<T_0=C_2\cdot\min \left(r^2, K^{-1} \right)$.
We will show that this leads to a contradiction.

First, notice that the curvature tensor $\Rm$ satisfies (\ref{9}),
then we have
\begin{align*}
\frac{\partial}{\partial
t}|\Rm|^2=\phi^2\triangle|\Rm|^2-2|\nabla\Rm|^2+\Rm*\Rm*\Rm.
\end{align*}
According to the proof of Theorem \ref{thm2}, we obtain
\begin{align*}
\sup_{x\in M}\|\Rm(g(t))\|_{2;B_r(x)} < 2\sup_{x\in
M}\|\Rm(g_0)\|_{2;B_r(x)}\le 2 [C_1A_0]^{-1},
\end{align*}
which implies a strict inequality for (\ref{14}).

Next, since the Ricci curvature satisfies (\ref{10}), then Corollary
\ref{cor11} and the covering argument imply that
\begin{align*}
|\Ric(g(t))|\le& CA_0^{\frac
{2}{3}}t^{-\frac{2}{3}}\left(\int_{B_{3r}(x)}|\Ric(g_0)|^3\right)^{\frac13}\\
                  \le&
                  Ct^{-\frac{2}{3}}\left(KA_0^2\int_{B_{3r}(x)}|\Ric(g_0)|^2\right)^{\frac{1}{3}}\\
                  \le& C_4t^{-\frac{2}{3}}
\end{align*}
where $\|\nabla\phi\|_{\infty}$ can be evaluated at $g(0)$, since
the metrics $g(t)$ are equivalent within the maximal time
$T_{\max}$.

Applying the bound on $\Ric$ to the following
$$\left|\frac{\der}{\der t}\int f^p\der V_g\right|\le 2\|\Ric\|_{\infty}\int f^p\der V_g,$$
we have
$$-2C_4t^{-\frac{2}{3}}\der t\le\der\log\int f^p\der V_g\le 2C_4t^{-\frac{2}{3}}\der t,$$
which implies that for some suitably chosen constants,
$$\left|\log\frac{\|f\|_p(t)}{\|f\|_p(0)}\right|<\log 2.$$
The differential inequality
$$\left|\frac{\der}{\der t}\int |\nabla f|^2\der V_g\right|\le 2\|\Ric\|_{\infty}\int |\nabla f|^2\der V_g$$
leads to a similar estimate.  Therefore, it follows that for any
$t\le T_0$,
$$\|f\|^2_4(t)<2\|f\|^2_4(0)\le 2A_0\|\nabla f\|^2_2(0)<4A_0\|\nabla f\|^2_2(t),$$
that is to say (\ref{12}) holds with strict inequality.

To show that (\ref{13}) holds with strict inequality, we use
Hamilton's trick.  Simply fix a tangent vector $v$ with respect to
$g(t)$, then
$$\frac{\der}{\der t}|v|^2_{g(t)}=\frac{\der}{\der t}(g_{ij}(t)v^i v^j)=g_{ij}^{\prime}(t)v^i v^j$$
implies
$$\left|\frac{\der}{\der t}\log |v|^2_{g(t)}\right|\le |g_{ij}^{\prime}(t)|\le 2|\Ric|.$$
So for $0\le t\le T_2 < T_0$,
$$\log \frac{|v|^2_{g(t)}}{|v|^2_{g(0)}}\le \int_0^{T_2}|g_{ij}^{\prime}(t)|\der t\le 2\|\Ric\|_{\infty}T_2<\log 2,$$
which implies
$$\frac{1}{2}|v|^2_{g(0)}<|v|^2_{g(t)}<2|v|^2_{g(0)},$$
for $t<T_0$.

Now we can show that $g(t)$ has a smooth limit as $t\rightarrow
T_{\max}$. If $T_{\max}<T_0$, we would be able to extend the
solution to (\ref 8) smoothly beyond $T_{\max}$ with (\ref{12}),
(\ref{13}) and (\ref{14}) still holding. This contradicts the
assumed maximality of $T_{\max}$. Hence, we conclude that
$T_{\max}\ge T_0$.

The estimate (\ref{11}) follows from Theorem \ref{thm2}. $\Box$

\vskip 0.2cm

Therefore, the previous theorem can be restated as the following.

\begin{thm}\label{thm5}
There exist constant $\varepsilon$ and $C_1$ such that for $r\le 1$,
if
$$\frac{r^4}{\Vol(B_r(x))}\int_{\B_r(x)}|\Rm(g_0)|^2\der V_{g_0}\le \varepsilon$$
and
$$|\Ric(g_0)|\le K,$$
then the equation (\ref{8}) has a smooth solution for $t\in[0,T)$,
where
$$T\ge C_1\cdot\min \left(r^2, K^{-1} \right).$$
Moreover, for $t\in(0,T)$, the Riemannian curvature tensor satisfies
the following bound,
\begin{align*}
\|\Rm\|_{\infty}\le C_2t^{-1}.
\end{align*}
Here $\varepsilon$, $C_1$ and $C_2$ only depend on the dimension of
$M$.
\end{thm}

\vskip 0.2cm

\noindent\textit{Proof: } When using the covering argument, we shall
use the following result. For any $y_1,\ y_2\in B_{\frac 32r}(x)$,
there exists a definite constant $c_1>0$ such that
$$\frac{1}{c_1}\Vol(B_r(y_1))\le\Vol(B_r(y_2))\le c_1\Vol(B_r(y_1)).$$
This follows from the volume comparison theorem, the facts that
within the maximal time all metrics are equivalent and that at $t=0$
$B_r(y_i)\cap B_r(x)\neq\emptyset$, $i=1,2$.

Then we can proceed as in the proof of the previous theorem. $\Box$

\vskip 0.2cm

Let $M$ be a 4-dimensional Riemannian manifold with bounded Ricci
curvature. Based on the result in \cite{ChGi}, we can show that if
there exists some $\varepsilon>0$ such that
$$\int_{B_r(x)}|\Rm|^2<\varepsilon,$$
then for some constant $c<1$, we have
$$\frac{r^4}{\Vol(B_{cr}(x))}\int_{B_{cr}(x)}|\Rm|^2<\varepsilon.$$
Therefore we obtain the main
theorem as required.

\end{document}